\documentclass[a4, 12pt]{amsart}
\usepackage[mathscr]{eucal}
\usepackage{amssymb}
\usepackage{latexsym}
\usepackage{amsthm}
\theoremstyle{plain}
\newtheorem{theorem}{Theorem}[section]

\newtheorem{lemma}{Lemma}[section]
\newtheorem{proposition}{Proposition}[section]

\makeatletter
\@addtoreset{equation}{section}

\title[2-dimensional complete self-shrinkers]
{2-dimensional complete self-shrinkers  in $\mathbf R^3$}
\author [Qing-Ming Cheng and Shiho Ogata]{Qing-Ming Cheng* and Shiho Ogata }
\dedicatory{Dedicated to Professor Yoshihiko Suyama for his 70th birthday}
\address{Qing-Ming Cheng \\  \newline \indent Department of Applied Mathematics, Faculty of Sciences, 
\newline \indent Fukuoka University, Fukuoka  814-0180, Japan.  \newline \indent cheng@fukuoka-u.ac.jp}
\address{Shiho Ogata \\  \newline \indent Department of Applied Mathematics, Graduate School of Sciences, 
\newline \indent Fukuoka University, Fukuoka  814-0180, Japan.  \newline \indent sd150501@cis.fukuoka-u.ac.jp}

\begin{document}
\maketitle

\begin{abstract}
It is our purpose to study complete self-shrinkers in Euclidean space.
First of all, we show some examples of complete self-shrinkers without
polynomial volume growth. By making use of the  generalized maximum
principle for $\mathcal{L}$-operator, we  give a complete classification
for 2-dimensional complete  self-shrinkers with constant squared norm of the 
second fundamental form in $\mathbb R^3$. In \cite{DX2},  Ding and Xin 
have proved this result under the assumption of  polynomial volume growth,
which is removed in our theorem.

\end{abstract}

\footnotetext{* Research partially Supported by JSPS Grant-in-Aid
for Scientific Research (B) No. 24340013 and Challenging Exploratory Research No. 25610016.}

\section{introduction}
\noindent
Let  $X: M^n\to \mathbb{R}^{n+1}$  be an $n$-dimensional hypersurface in the $n+1$-dimensional
Euclidean space $\mathbb{R}^{n+1}$. If the position vector $X$ evolves in the direction of the mean
curvature $H$, then it gives rise  to a solution to mean curvature flow:
$$
X(\cdot, t):M^n\to  \mathbb{R}^{n+1}
$$
satisfying $X(\cdot, 0)=X(\cdot)$ and
\begin{equation}
\dfrac{\partial X(p,t)}{\partial t}=H(p,t), \quad (p,t)\in M\times [0,T),
\end{equation}
where $H(p,t)$ denotes the mean curvature vector  of hypersurface $M_t=X(M^n,t)$ at point $X(p,t)$.
The equation (1.1) is called the mean curvature flow equation.
The study of the mean curvature  flow from the perspective of partial differential
equations commenced with Huisken's paper  \cite{H1}  on the flow of convex
hypersurfaces (cf. \cite{EH}). 

\noindent
One of the most important problems in the  mean curvature flow is to understand
the possible singularities that the flow goes through.  A key starting point
for singularity analysis is Huisken's monotonicity formula because the monotonicity
implies that the flow is asymptotically self-similar near a given singularity and
thus, is modeled  by  self-shrinking solutions of the flow.

\noindent
An  $n$-dimensional  hypersurface  $X: M\rightarrow \mathbb{R}^{n+1}$  in the $(n+1)$-dimensional
Euclidean space $\mathbb{R}^{n+1}$  is called a self-shrinker if it satisfies 
\begin{equation*}
H+\langle X, N\rangle=0,
\end{equation*}
where  $H$ and $N$ denote  the mean curvature  and the unit normal vector of the hypersurface, respectively.
It is known that self-shrinkers play an important role in the study of the mean curvature flow because
they describe all possible blow up at a given singularity of the  mean curvature flow. For classifications of complete
self-shrinkers, Abresch and Langer \cite{AL},  Huisken \cite{H2, H3} and  Colding and Minicozzi \cite{CM1} have 
obtained  very important results. In fact,  Abresch and Langer \cite{AL}  classified 
closed self-shrinker curves in $\mathbb{R}^2$ and showed that the  round circle is the only embedded self-shrinkers.
Huisken \cite{H2, H3} and  Colding and Minicozzi \cite{CM1}  have proved that  if $X: M\rightarrow \mathbb{R}^{n+1}$ 
 is  an $n$-dimensional  complete  embedded self-shrinker
in $\mathbb{R}^{n+1}$ with $H\geq 0$ and  with polynomial volume growth,
then $X: M\rightarrow \mathbb{R}^{n+1}$  is  isometric to either $\mathbb{R}^{n}$,  the  round sphere $S^{n}(\sqrt{n})$, or a cylinder $S^m (\sqrt{m})\times \mathbb{R}^{n-m}$, $1\leq m\leq n-1$.
In \cite{CL}, Cao  has conjectured   that complete self-shrinkers must have polynomial volume growth. Furthermore, Ding and Xin \cite{DX1} and
X. Cheng and Zhou \cite{CZ} have proved that a  complete self-shrinker has  polynomial volume growth if and only if it is proper. 
From the following proposition, we know that there are many complete self-shrinkers without polynomial volume growth.
\begin{proposition}  
For any integer $n>0$, there exist $n$-dimensional  complete self-shrinkers without polynomial volume growth  in $\mathbb R^{n+1}$.
\end{proposition}

\noindent
In fact, in \cite{H}, Halldorsson has proved there exist complete self-shrinker curves $\Gamma$ in $\mathbb R^2$, which is contained in an annulus around the origin 
and whose image is dense in the  annulus. Hence, these complete self-shrinker curves $\Gamma$  are not proper. Thus, for any integer $n>0$,
$\Gamma \times \mathbb R^{n-1}$ is  a complete self-shrinker in $\mathbb R^{n+1}$, which  does not have polynomial volume growth.

\vskip1mm
\noindent
In \cite{CL}, Cao and Li have proved that if an $n$-dimensional  complete  self-shrinker $X: M\rightarrow \mathbb{R}^{n+1}$ with  polynomial volume growth satisfies $S\leq 1$, then   $X: M\rightarrow \mathbb{R}^{n+1}$  is  isometric to either  $\mathbb{R}^{n}$,  the  round sphere $S^{n}(\sqrt{n})$, or a cylinder $S^m (\sqrt{m})\times \mathbb{R}^{n-m}$, $1\leq m\leq n-1$ (cf. \cite{LW}).
Furthermore, Ding and Xin \cite{DX2} have studied 2-dimensional complete self-shrinkers with polynomial volume growth and with constant  squared norm of the second fundamental form (cf. \cite{DX2} and \cite{CW} for any dimension). They have proved

\noindent
{\bf Theorem DX}. {\it 
Let $X: M\rightarrow \mathbb{R}^{3}$  be a $2$-dimensional complete self-shrinker  with polynomial volume growth
in $\mathbb{R}^{3}$.
If the squared norm $S$ of the second fundamental form is constant,
 then $X: M\rightarrow \mathbb{R}^{3}$  is isometric to one of the following:

\begin{enumerate}
\item $\mathbb{R}^{2}$,
\item
a cylinder $S^1 (1)\times \mathbb{R}$
\item the round sphere $S^{2}(\sqrt{2})$.
\end{enumerate}
}

\noindent
In this paper, we want to remove the assumption of polynomial volume growth in the above theorem of Ding and Xin 
and to prove that  the above result  of Ding and Xin holds by making use of a different method.

\begin{theorem}
Let $X: M\rightarrow \mathbb{R}^{3}$  be a $2$-dimensional complete self-shrinker in $\mathbb{R}^{3}$.
If the squared norm $S$ of the second fundamental form is constant,
 then $X: M\rightarrow \mathbb{R}^{3}$  is isometric to one of the following:
\begin{enumerate}
\item $\mathbb{R}^{2}$,
\item
a cylinder $S^1 (1)\times \mathbb{R}$
\item the round sphere $S^{2}(\sqrt{2})$.
\end{enumerate}\end{theorem}

\section{Proof of theorem 1.1}

\noindent 
Let $X: M^2\rightarrow\mathbb{R}^{3}$ be a 2-dimensional surface in 
$\mathbb{R}^{3}$. We choose a local orthonormal frame field
$\{e_A\}_{A=1}^{3}$ in $\mathbb{R}^{3}$ with dual co-frame field
$\{\omega_A\}_{A=1}^{3}$, such that, restricted to $M^2$,
$e_1, e_2$ are tangent to $M^2$. Hence,  we have
\begin{equation*}
dX=\sum_{i=1}^2\omega_i e_i, \quad de_i=\sum_{j=1}^2\omega_{ij}e_j+\omega_{i3}e_3.
\end{equation*}
We restrict these forms to $M^2$, then
\begin{equation}
\omega_3=0 
\end{equation}
and 
\begin{equation*}
\omega_{i3}=\sum_{j=1}^2h_{ij}\omega_j,\quad
h_{ij}=h_{ji},
\end{equation*}
where $h_{ij}$ denote components of the second fundamental form of $X: M^2\rightarrow\mathbb{R}^{3}$.
Take $e_1, e_2$ such that, at any fixed point, 
 $$
 h_{ij}=\lambda_i\delta_{ij},
 $$
where $\lambda_1$ and $\lambda_2$ are the principal curvatures of $X: M^2\rightarrow\mathbb{R}^{3}$.
Thus, the Gauss curvature $K$ and the mean curvature $H$  are given by
$$
K=\lambda_1\lambda_2, \quad H=\lambda_1+\lambda_2.
$$
For a smooth function $f$, the $\mathcal{L}$-operator is defined by
\begin{equation}
\mathcal{L}f=\Delta f-\langle X,\nabla f\rangle
\end{equation}
where $\Delta$ and $\nabla$ denote the Laplacian and the gradient
operator on the self-shrinker, respectively.
In order to  prove  our results,  the following
generalized maximum principle for $\mathcal{L}$-operator on self-shrinkers is very important, which is proved
by Cheng and Peng in \cite{CP}:
\begin{lemma} {\rm(}Generalized maximum principle for $\mathcal{L}$-operator {\rm)}
Let $X: M^n\to \mathbb{R}^{n+p}$ {\rm (}$p\geq 1${\rm)} be a complete self-shrinker with Ricci
curvature bounded from below. Let $f$ be any $C^2$-function bounded
from above on this self-shrinker. Then, there exists a sequence of points
$\{p_k\}\subset M^n$, such that
\begin{equation*}
\lim_{k\rightarrow\infty} f(X(p_k))=\sup f,\quad
\lim_{k\rightarrow\infty} |\nabla f|(X(p_k))=0,\quad
\limsup_{k\rightarrow\infty}\mathcal{L} f(X(p_k))\leq 0.
\end{equation*}
\end{lemma}

\noindent
\emph{Proof of theorem 1.1}.
Since $X: M^2\rightarrow\mathbb{R}^{3}$  is a complete self-shrinker, we have
\begin{equation}
H+\langle X, N \rangle=0.
\end{equation}
By a simple calculation, we have
\begin{equation*}
\dfrac12\mathcal LS=\sum_{i,j,k}h_{ijk}^2+S(1-S),
\end{equation*}
where $S=\sum_{i,j=1}^2h_{ij}^2$ is the squared norm of the second fundamental form
and $h_{ijk}$ denote components of the first covariant derivative of the second fundamental form.
Since $S$ is constant, we have
\begin{equation}
\sum_{i,j,k}h_{ijk}^2+S(1-S)=0.
\end{equation}
If $S=1$, then we know $h_{ijk}\equiv0$. Hence, 
$X: M^2\rightarrow\mathbb{R}^{3}$ is isometric to the round sphere $S^2(\sqrt{2})$
or the cylinder $S^1(1)\times \mathbb{R}$ from the results of Lawson \cite{L}.
If $S<1$, from the theorem of Cheng and Peng \cite{CP}, we know that 
$X: M^2\rightarrow\mathbb{R}^{3}$ is isometric to $\mathbb R^2$.

\noindent
Next, we  prove that $S\leq 1$ holds. 
By a direct computation,  we have
\begin{equation}
\dfrac12\mathcal L|X|^2=2-|X|^2.
\end{equation}
Since $S$ is constant, we know that the Gauss curvature satisfies
$$
K=\lambda_1\lambda_2\geq -\frac{\lambda_1^2+\lambda_2^2}2=-\frac{S}2.
$$
Therefore, the Gauss curvature is bounded from below.
Since $-|X|^2\leq 0$ is bounded from above, 
we can apply the generalized maximum principle for $\mathcal L$-operator
to the function $-|X|^2$. Thus, there exists a sequence $\{p_k\}$ in $M^2$ such that 
\begin{equation}
\lim_{k\rightarrow\infty} |X|^2(p_k)=\inf |X|^2,\quad
\lim_{k\rightarrow\infty} |\nabla |X|^2(p_k)|=0,\quad
\liminf_{k\rightarrow\infty}\mathcal{L} |X|^2(p_k)\geq 0.
\end{equation}
From (2.5) and (2.6), we have
\begin{equation}
\inf |X|^2\leq 2.
\end{equation}
Since $|\nabla |X|^2|=\sum_{i=1}^2\langle X, e_i\rangle^2$ holds, we have from (2.6)
\begin{equation*}
\lim_{k\rightarrow\infty} |\nabla |X|^2(p_k)|=\lim_{k\rightarrow\infty}\sum_{i=1}^2\langle X, e_i\rangle^2(p_k)=0.
\end{equation*}
Hence, we get from (2.3)
\begin{equation}
\inf|X|^2=\lim_{k\rightarrow\infty}H^2(p_k), \quad \lim_{k\rightarrow\infty}|\nabla H|(p_k)=0.
\end{equation}
Since $S$ is constant, from the definition of the mean curvature $H$ and (2.3), we obtain, for $j=1, 2$,
\begin{equation}
\begin{aligned}
&\lim_{k\rightarrow\infty}\big(h_{11j}(p_k)+h_{22j}(p_k)\big)=0,\\
&\lim_{k\rightarrow\infty}\big(\lambda_1(p_k)h_{11j}(p_k)+\lambda_2(p_k)h_{22j}(p_k)\big)=0.
\end{aligned}
\end{equation}
Since $S$ is constant, we know that $\{\lambda_j(p_k)\}$ and $\{h_{iij}(p_k)\}$  are bounded sequences.
Thus, we can assume
$$
\lim_{k\rightarrow\infty}h_{iij}(p_k)=\bar h_{iij}, \quad \lim_{k\rightarrow\infty}\lambda_j(p_k)=\bar \lambda_j,
$$
for $i, j=1, 2$.
From (2.9), we obtain 
\begin{equation}
\begin{cases}
 \bar h_{11j}+\bar h_{22j}&=0,\\
\bar \lambda_1\bar h_{11j}+\bar \lambda_2\bar h_{22j}&=0.
\end{cases}
\end{equation}
If $\bar \lambda_1\neq \bar \lambda_2$ is satisfies, according to (2.10), we infer
$$
\bar h_{iij}=0
$$
for $i, j=1, 2$. According to Codazzi equations, we have 
\begin{equation*}
\sum_{i,j,k}\bar h_{ijk}^2=0.
\end{equation*}
From (2.4), we have $S=1$ or $S=0$. Hence  $S\leq 1$.

\noindent
If $\bar \lambda_1= \bar \lambda_2$ holds, we have
$$
S=\bar\lambda_1^2+\bar\lambda_2^2=\frac{\bigl(\bar\lambda_1
+\bar\lambda_2\bigl)^2}2=\frac{\lim_{k\rightarrow\infty}H^2(p_k)}2.
$$
According to (2.7) and (2.8), we have
$$
S\leq 1.
$$
Hence, $S=0$ or $S=1$.
According  to the theorem of Lawson \cite{L},  we know that
$M^n$ is isometric to the round sphere $S^2(\sqrt{2})$, the cylinder $S^1(1)\times \mathbb{R}^{1}$ or $\mathbb R^2$.
\begin{flushright}
$\square$
\end{flushright}
\vskip 5mm

\vskip3mm
\noindent
{\bf Acknowledgement}. Authors would like to thank professor Wei Guoxin for fruitful discussions.

\end{document}